# The persistence exponents of Gaussian random fields connected by the Lamperti transform

## G. Molchan

[1] Institute of Earthquake Prediction Theory and Mathematical Geophysics,
Russian Academy of Science, 84/32 Profsoyuznaya st., Moscow, RF

[2] HSE University, Myasnitskaya str. 20, Moscow, RF

E-mail address: molchan@mitp.ru

**Abstract**. The (fractional) Brownian sheet is a simplest example of a Gaussian random field $X$ whose covariance is the tensor product of a finite number (d) of nonnegative correlation functions of self- similar Gaussian processes. Let $Y$ be the homogeneous Gaussian field obtained by applying to X the Lamperti transform, which involves the exponential change of time and the amplitude normalization to have unit variance. Under some assumptions, we prove the existence of the persistence exponents for both fields, $X$ and $Y$, and find the relation between them. The exponent for $X$ is the limit of $-\ln P(X(\mathbf{t}) \leq 1, \mathbf{t} \in [0,T]^d)/\ln^d T$, $T \gg 1$. In terms of $Y$ it has the form $-\lim \ln P(Y(\mathbf{t}) \leq 0, \mathbf{t} \in TG)/T^d := \theta$, $T \gg 1$, where $G$ is a suitable d-simplex and $T$ is a similarity coefficient; $G$ can be selected in form $[0,c]^d$ if d=2. The exponent $\theta$ exists for any continuous Gaussian homogeneous field $Y$ with non-negative covariance when $G$ is bounded region with a regular boundary. The problem of the existence and connection of exponents for $X$ and $Y$ was raised by Li and Shao [Ann. Probab. 32:1, 2004] and originally concerned the Brownian sheet.

## 1. Introduction

The *tensor product* of centered Gaussian processes $X_i(t)$ with covariances $B_i$ is by definition a Gaussian field $X(\mathbf{t})$, $\mathbf{t} = (t_1,...,t_d)$ with the covariance $B_X(\mathbf{t},\mathbf{s}) = \prod B_i(t_i, s_i)$. The simplest example, which is important in applications, is the Brownian sheet $W(\mathbf{t})$ or the Chentsov field, for which $B_i(t,s) = t \wedge s$. Csa'ki et al [6] drew attention to the unusual log-asymptotics of probability that $W(\mathbf{t})$ is bounded in $[0,T]^2$ for large $T$. The exact order of the asymptotics, $\psi_T$, was found by Li and Shao, [9]. They considered a more general case where $X_i(t), i=1-d$ are fractional Brownian motions $W_h(t)$, $0 < h < 1$ with

$$B_i(t,s) = (t^{2h} + s^{2h} - |t-s|^{2h})/2 \ , \ i = 1-d. \tag{1.1}$$

In this case, it was shown that $\psi_T = \ln^d T$ and

$$0 < \theta_X^-(1+o(1)) \leq -\ln P(X(\mathbf{t}) \leq 1, \mathbf{t} \in [0,T]^d)/\psi_T \leq \theta_X^+(1+o(1)), \quad T \to \infty,$$

where $\theta_X^\pm$ are unknown finite constants.

Li and Shao raised two questions. The first is the existence of a persistence exponent for $X(\mathbf{t}) = \prod \otimes W_h(t_i)$, i.e., the equality $\theta_X^- = \theta_X^+$. The positive solution for d=2 was given in [12]. The second question is related to the persistence exponents of Gaussian random fields $X$ and $Y$ connected by the *Lamperti transform*:

$$Y(\mathbf{t}) = X(e^{\mathbf{t}})/[EX^2(e^{\mathbf{t}})]^{1/2}, \quad e^{\mathbf{t}} = (e^{t_1},...,e^{t_d}). \tag{1.2}$$

Consider for a moment one-dimensional time. If $\{X(t), 0 \leq t \leq T, EX^2(1) = 1\}$ is a self-similar process of index h (h-ss), i.e. $X(\lambda t) =_{law} \lambda^h X(t)$, then $Y(t) = X(e^t)e^{-th}, t \leq \ln T$ and $Y(t)$ is stationary. Since $EX^2(\mathbf{t}) = \prod_{i=1}^d EX_i^2(t_i)$, the same is true for the tensor product of $h_i - ss$ processes $X_{h_i}$, that is, $Y(\mathbf{t})$ is homogeneous; moreover, $Y(\mathbf{t})$ is the tensor product of $Y_{h_i}(\cdot)$ related to $X_{h_i}$.

Reducing the problem for a self-similar object $X$ to a homogeneous one is quite often used in applications, [8, 9,11,12,14]; for example, for a large class of h-ss Gaussian processes the dual pair $X_h$ and $Y_h$ has the following identical limits:

$$\lim_{T\to\infty} -\ln P(X_h(t) \leq 1, \mathrm{t} \in [0,T])/\ln T$$
$$= \lim_{T\to\infty} -\ln P(Y_h(t) \leq 0, t \in [0,T])/T, \tag{1.3}$$

The second question by Li and Shao [9] is to find an analogue of (1.3) for the Brownian sheet.

Below, we consider Gaussian fields of the tensor product type with arbitrary ss- components to give an answer to Li-Shao's questions in the general situation. Along the way, we generalize the fact of the existence of the exponent (1.3) for any continuous Gaussian stationary processes with non-negative covariance, [9], to homogeneous fields.

## 2. Results.

Below we will use the following notation and properties:

$X(\mathbf{t}) = \prod_{i=1}^d \otimes X_i(t_i)$ is a centered Gaussian field represented by the tensor product of centered Gaussian processes $X_i, i = 1-d$, i.e. $EX(\mathbf{t})X(\mathbf{s}) = \prod_{i=1}^d EX_i(t_i)X(s_i)$;

if $X_i$, $EX_i^2(1) = 1$ components are self-similar with indexes $h_i$ then for any fixed $\{\lambda_i > 0, i = 1-d\}$

$$X(\lambda_1 t_1,...,\lambda_d t_d) =_{law} \lambda_1^{h_1}...\lambda_d^{h_d} X(\mathbf{t}), \tag{2.1}$$

$$X(\lambda_1,...\lambda_{k-1}, s\lambda_k, \lambda_{k+1}...,\lambda_d) =_{law} \lambda_1^{h_1}...\lambda_d^{h_d} X_k(s). \tag{2.2}$$

These relations follow from the equality of covariances of the compared random functions.

Further, $M_X(G) = \sup\{X(\mathbf{t}), \mathbf{t} \in G\}$; $G \in R^d$ is a region (connected set having internal points ) of finite volume $|G|$; $TG$ is a region similar to G with a similarity coefficient $T$.
The limit (if it exists)

$$\theta_X[G, \psi_T, c] = \lim_{T \to \infty} -\ln P(M_X(TG) \le c)/\psi(T) \qquad (2.3)$$

is the *persistence exponent* of $X(\mathbf{t})$ relative to $TG$. To simplify the notation, we will use the symbols $\theta_X(G)$ for the exponent $\theta_X[G, \ln^d T, 1]$ if $X(\mathbf{t})$ a self-similar field and $\theta_Y(G)$ for the exponent $\theta_Y[G, T^d, 0]$ if $Y(\mathbf{t})$ is a homogeneous random field. Similarly, the symbols $\theta_{\bullet}^-(G), \theta_{\bullet}^+(G)$ will correspond to the operations $\liminf$ and $\limsup$ in (2.3) respectively.

### 2.1 Persistence exponents for homogeneous fields.

We will consider bounded regions $G \subset R^d$ with a *regular* boundary. This means that on the standard rectangular lattice of step $\delta$, the internal and external measures of $G$, associated with the summation of the lattice elements, are equal to $|G|(1+o(1))$ as $\delta \downarrow 0$.

**Theorem 1.** Let $Y(\mathbf{t}), \mathbf{t} \in R^d$ be a continuous, centered, Gaussian homogeneous field with nonnegative covariance, $B_Y(\mathbf{t}) \ge 0$, $B_Y(\mathbf{0}) = 1$. Then, for any bounded region $G$ with a *regular* boundary, the persistence exponent, $\theta_Y(G)$, exists and is finite. If, in addition, $B_Y(\mathbf{t}) \in L_1$, then $\theta_Y(G) > 0$.

**Remark.** A similar result for a stationary process with $B_Y(t) \ge 0$ was first obtained in [9] as a simple consequence of the subadditivity of the function: $T \to -\ln P(Y(t) < c, t \in (0,T))$. This property is lost when $G$ is an arbitrary connected region in $R^d, d > 1$.

The proof of the theorem will be given after the formulation of auxiliary statements of independent interest.

As follows from the definition of the exponent $\theta_Y(G)$ for homogeneous fields, we have

$$\theta_Y(kG) = k^d \theta_Y(G), \qquad (2.4)$$

provided that one of the indexes exists. Let's consider a less trivial case of the relationship of exponents for different regions. Assume that for large $T$, the region $TG_1$ can be approximated from above by combining the regions $\{\gamma_i G_2, i = 1 - N_T\}$, i.e., $TG_1 \subseteq \cup_1^{N_T} \gamma_i G_2$; here $\{\gamma_i\}$ are the elements of some group $\Gamma$ of movements of $R^d$. If $|TG_1| = N_T |G_2|(1+o(1))$, $T \to \infty$, then we will say that $G_2$ is a *parquet tile* for the figure $G_1$ relative to $\Gamma$.

**Proposition 2.** Let $Y(\mathbf{t}), \mathbf{t} \in R^d$ be a centered, Gaussian field with nonnegative covariance. Let $Y(\mathbf{t})$ be homogeneous with respect to some movements group $\Gamma$ of $R^d$, i.e. $Y(\gamma \mathbf{t}) =_{law} Y(\mathbf{t})$ for any fixed $\gamma \in \Gamma$. If the bounded regions $G_1$ and $G_2$ are parquet tiles for each other relative to $\Gamma$ then the exponents $\theta_Y(G_i), i = 1,2$ exist and

$$\theta_Y(G_1)/|G_1| = \theta_Y(G_2)/|G_2| \qquad (2.5)$$

(for the proof, see section 3).

**Corollary 3.** The d-cube $K = [0,1]^d$ and d-parallelepiped $\Pi = [0, a_1] \times ... \times [0, a_d]$ are a mutual parquet pair relative to the group of shifts of $R^d$. Therefore, the exponents $\theta_Y(K)$ and $\theta_Y(\Pi)$ for a homogeneous Gaussian field $Y(\mathbf{t}), \mathbf{t} \in R^d$, $B_Y(\mathbf{t}) \geq 0$ exist and $\theta_Y(\Pi) = \theta_Y(K) \prod_1^d a_i$.

**Proposition 4,** ([7]). Let $Y(\mathbf{t}), t \in R^d$ be a continuous, centered, Gaussian homogeneous field with a spectral measure having a density $f$ in some $\varepsilon$–vicinity of zero frequency, where $0 < c \leq f \leq C < \infty$. Then $\theta_Y^-(K) > 0$ and $\theta_Y^+(K) < \infty$.

**Remark.** The paper [7] refers to the case $d = 1$, although the proof method is easily adapted to the general case: $d \geq 1$.

**Proof of Theorem 1.**

Here and further we will use the notation: $K = [0,1]^d$, $G^c = K \setminus G$, and $P_Y(G) = P(M_Y(G) \leq 0)$. Consider a bounded region $G \subset R^d$ with a regular boundary. Thanks to (2.4), we can assume that $G \subset K$. Since the covariance $B_Y$ of the random field $Y(\mathbf{t})$ is nonnegative, Slepian's lemma, [10], implies

$$P_Y(TK) \geq P_Y(TG) P_Y(TG^c),$$
$$-\ln P_Y(TK)/T^d \leq -\ln P_Y(TG)/T^d - \ln P_Y(TG^c)/T^d. \qquad (2.6)$$

According to Corollary 1, the left part of (2.6) has the limit $\theta_Y(K)$ as $T \to \infty$. We choose a sequence $T \to \infty$ such that $-\ln P_Y(TG)/T^d \to \theta_Y^-(G)$. In that case (2.6) gives

$$\theta_Y(K) \leq \theta_Y^-(G) + \theta_Y^+(G^c) \qquad (2.7)$$

Let $\hat{G}_T \supset TG$ is the minimal approximation of $TG$ from above on rectangular lattice in $R^d$ with the unit step and $N_T = |\hat{G}_T|$. Given Slepian's lemma and the homogeneity of $Y(\mathbf{t})$, we have for any $c > 0$

$$P_Y(T(cG)) \geq P_Y(c\hat{G}_T) \geq [P_Y(cK)]^{N_T}.$$

The boundary regularity of $G$ gives: $N_T = |\hat{G}_T| = |TG|(1+o(1))$. Therefore

$$-\ln P_Y(TcG)/(Tc)^d \leq -|G|[P_Y(cK)]/c^d (1+o(1)), \quad o(1) \to 0 \text{ as } T \to \infty.$$

Passing to infinity, first over $T$, and then over $c$, we get

$$\theta_Y^+(G) \leq |G|\theta_Y(K).\qquad(2.8)$$

The boundary of $G^c = K \setminus G$ is regular, since it includes the boundaries of $K$ and $G$. Therefore, similarly, we have

$$\theta_Y^+(G^c) \leq |G^c|\theta_Y(K).\qquad(2.9)$$

Since $|G| + |G^c| = |K| = 1$,

$$\theta_Y^+(G) + \theta_Y^+(G^c) \leq \theta_Y(K).\qquad(2.10)$$

Combining (2.7) and (2.10), we get

$$\theta_Y(K) \leq \theta_Y^-(G) + \theta_Y^+(G^c) \leq \theta_Y^+(G) + \theta_Y^+(G^c) \leq \theta_Y(K).$$

Therefore

$$\theta_Y^-(G) + \theta_Y^+(G^c) = \theta_Y^+(G) + \theta_Y^+(G^c) = \theta_Y(K).$$

Assuming that $\theta_Y(K) < \infty$, we get $\theta_Y^-(G) = \theta_Y^+(G) < \infty$.

We now are going to show that the assumption is valid

As above, we have

$$-\ln P_Y(N\varepsilon K)/N^d \leq -\ln P_Y(\varepsilon K).$$

This gives $\theta_Y(\varepsilon K) \leq -\ln P_Y(\varepsilon K)$ when $N \to \infty$. Assuming $\theta_Y(K) = \infty$, we have $P_Y(\varepsilon K) = 0$. Since $P_Y(\varepsilon K) = P_{-Y}(\varepsilon K)$, the continuity of $Y(\cdot)$ implies $P\{\exists \mathbf{t} \in \varepsilon K : Y(\mathbf{t}) = 0\} = 1$ for any $\varepsilon > 0$. In that case, however, we have $Y(\mathbf{0}) = 0$ a.s., which contradicts the condition $EY^2(\mathbf{0}) = 1$.

The non-triviality of $\theta_Y(G)$, i.e. $\theta_Y(G) > 0$, requires additional assumptions. Let $B_Y(\mathbf{t}) \in L_1$. Then $Y(\mathbf{t})$ has a continuous spectral function

$$f(\lambda) = (2\pi)^{-d}\int \cos(\langle \lambda, \mathbf{t}\rangle)B_Y(\mathbf{t})d\mathbf{t} \geq 0, \lambda \in R^d.$$

Since $B_Y(\mathbf{t}) \geq 0$, $f(\lambda) \leq (2\pi)^{-d}\|B_Y\|_{L_1} = f(0) < \infty$. Because of the continuity of $f(\lambda)$, there is a neighborhood of $\lambda = 0$, where $0 < c < f(\lambda) \leq f(0)$. Hence, by Proposition 2, $\theta_Y^-(K) > 0$ and $\theta_Y^+(K) < \infty$. For $G \supseteq \varepsilon K, P_Y(TG) \leq P_Y(T\varepsilon K)$ and therefore $\theta_Y(G) \geq \theta_Y(\varepsilon K) = \varepsilon^d \theta_Y(K) > 0$.

**2.2 Persistence exponents of self-similar fields.**

Let $X_i$ be a self similar Gaussian process, and $Y_i$ be its stationary double obtained by the Lamperti transform. This section discusses the Gaussian fields, symbolically represented by the tensor products $X(\mathbf{t}) = \prod_{i=1}^{d} \otimes X_i(t_i)$ and $Y(\mathbf{t}) = \prod_{i=1}^{d} \otimes Y_i(t_i)$. We will make the following assumptions about the components $\{X_i, i = 1-d\}$:

a) $X_i$ is a $h_i$- ss centered Gaussian process, $EX_i^2(1) = 1$, and is time-invertible, that is, $X_i(t) =_{law} X_i(1/t)t^{2h_i}, t > 0$, $h_i > 0$;

b) $E|X_i(t) - X_i(s)|^2 \le c|t-s|^{2\alpha_i}$ for $t,s \subset [0,1]$, and $0 < \alpha_i \le 1$;

c) the covariance of $X_i$ is nonnegative: $B_i(t,s) = EX_i(t)X_i(s) \ge 0$;

d) $B_i(1,s) \le s^{h_i} / \ln^{1+\varepsilon} s$, $s > s_0 > 0$ for some $\varepsilon > 0$;

e) $B_i(1, 1+s) \ge c_i > 0$, $s > 0$.

*Remarks*. 1. Condition (a) means that the dual process $Y_i$ is time reversible: $Y_i(t) =_{law} Y_i(-t)$;

2. Condition (b) guarantees the continuity (a.s.) of $X_i$ on $[0,1]$. This is also true for any interval $[0,a]$ because $X_i(at) =_{law} a^{h_i} X_i(t)$;

3. Under the conditions (d,c), the dual process $Y_i$ has a continuous spectral density;

4. The fractional Brownian motion $W_h(t)$, $0 < h < 1$, for which $E|W_h(t) - W_h(s)|^2 = |t-s|^{2h_i}$, satisfies all the conditions (a-e) with constants $\alpha = h, c = 1, \varepsilon > 0$.

*Notation:* $\{\mathbf{t} \ge \mathbf{a}\} = \{\mathbf{t} = (t_1,...,t_d): t_i \ge a_i, i = 1-d\}$, $\mathbf{h} = (h_1,...,h_d)$, $H = \sum_1^d h_i$.

**Theorem 5**. Let $K = [0,1]^d$ and $S_\mathbf{h} = \{\mathbf{t} \in R^d : \mathbf{t} \ge \mathbf{0}, (\mathbf{t},\mathbf{h}) \le H\}$. Consider the pair of dual Gaussian fields $X(\mathbf{t}), \mathbf{t} \in R_+^d$ and $Y(\mathbf{t}), \mathbf{t} \in R^d$ under conditions (a-e). Then the persistence exponents $\theta_X(K)$ and $\theta_Y(S_\mathbf{h})$ exist, are equal, and are nontrivial (nonzero and finite). In the case of $d = 2$ we also have

$$\theta_X(K) = \theta_Y(K) \cdot 2(h_a/h_g)^2 = \theta_Y(cK), \tag{2.11}$$

where $h_a = (h_1 + h_2)/2$, $h_g = \sqrt{h_1 h_2}$, and $c = \sqrt{2} h_a / h_g$.

This statement is based on Propositions 6, 7. Their proofs are postponed until Section 3.

**Proposition 6.** Let $U_T = \{\mathbf{t} = (t_1,...t_d): t_1^{h_1} \times ... \times t_d^{h_d} \le 1\} \cap [0,T]^d$. Then under conditions (a,b,c), we have

$$EM_X(U_T) \le c\sqrt{\ln T}. \tag{2.12}$$

**Proposition 7.** Let $(\mathbf{H}_T, \|\cdot\|_{\mathbf{H}_T})$ be a Hilbert space of functions on $TK \subset R^d$ with the reproducing kernel $B_X(\mathbf{t},\mathbf{s}) = EX(\mathbf{t})X(\mathbf{s})$, related to $X(\mathbf{t}) = \prod_{i=1}^d \otimes X_i(t_i)$, [2]. Then under the conditions (a-e) there exists $\varphi(\mathbf{t}) \in \mathbf{H}_T$ such that

$$\varphi(\mathbf{t}) \ge 1, \mathbf{t} \in U_T^c \text{ and } \|\varphi\|_{\mathbf{H}_T} \le c \ln^{(d-1)/2} T. \tag{2.14}$$

where

$$U_T^c = \{\mathbf{t} = (t_1,...t_d): t_1^{h_1} \times ... \times t_d^{h_d} > 1\} \cap [0,T]^d. \tag{2.15}$$

**Proof of Theorem 5**.

Our goal is to show that the asymptotics of $P_X(TK), T \to \infty$ is determined by the asymptotics of $P(X(\mathbf{t}) \leq 0, \mathbf{t} \in U_T^c\}$. The latter is reduced by the Lamperti transform to the asymptotics of $P_Y(\tilde{T}S_\mathbf{h}), \tilde{T} = \ln T \to \infty$ for the homogeneous field $Y$; this is the case when we can apply Theorem 1.

*Step 1.* We show that under conditions (a,b,c), $\theta_X^+(K) \leq \theta_X^+[U_T^c, \ln^d T, 0]$.

For any $C_T > 0$,

$$P(M_X(U_T^c) \leq 0) \leq P(M_X([0,T]^d) \leq C_T) + P(M_X(U_T) \geq C_T). \tag{2.16}$$

Applying the concentration principle to $M_X(U_T)$ (see e.g. [10, Ch.14]) with $C_T^2 = 4A \ln^d T$, we get for $A > 0$ and large T

$$R_T := P(M_X(U_T) \geq C_T) \leq \exp[-(C_T - EM_X(U_T))^2/(2\sigma^2)],$$

where $\sigma^2 = \max_{U_T} E[X(\mathbf{t})]^2 = 1$ and $EM_X(U_T) \leq c\sqrt{\ln T}$ (Proposition 6).

As a result

$$R_T \leq \exp(-A \ln^d T), T > T_0.$$

Because of the self-similarity of $X(\mathbf{t})$,

$$P(M_X([0,T]^d) \leq C_T) = P(M_X([0, \lambda T]^d) \leq 1),$$

where $\lambda^H C_T = 1$ or $\lambda = C_T^{-1/H} = (4A\ln T)^{-1/H}$.

Hence, if $\theta_+ := \theta_X^+(V_T, \ln^d T, 0) < \infty$ and $A \geq \theta_+$, then $R_T = o(P(M_X(U_T^c) \leq 0))$.

Therefore

$$-\frac{\ln P(M_X(U_T^c) \leq 0) + o(1)}{\ln^d T} \geq -\frac{\ln P(M_X([0, \lambda T]^d) \leq 1)}{\ln^d \lambda T} \frac{\ln^d \lambda T}{\ln^d T}.$$

Since $\ln \lambda = o(\ln T)$, we get

$$\theta_X^+[U_T^c, \ln^d T, 0] \geq \theta_X^+(K). \tag{2.17}$$

It remains to make sure that the left side of this inequality is finite.

Since the amplitude normalization of the process X does not affect the probability of exceeding the zero level, we have

$$p_T := P(M_X(U_T^c) \leq 0) = P(M_Y(V_{\tilde{T}}) \leq 0), \tilde{T} = \ln T,$$

where $V_{\tilde{T}} = \{\mathbf{t} : \mathbf{t} \leq \tilde{T}\mathbf{e}, (\mathbf{h},\mathbf{t}) \geq 0\}, \mathbf{e} = (1,...,1)$.

The (a)- assumption means that $Y$ is invariant with respect to time shifts and axial-symmetries of the form: $\mathbf{t} \to (\varepsilon_1 t_1,..., \varepsilon_d t_d), \varepsilon_i = \pm 1$. The described group of movements $\Gamma$ allows us to convert $V_{\tilde{T}}$ to $\tilde{T}S_\mathbf{h}$ where $S_\mathbf{h} = \{\mathbf{t} \in R^d : \mathbf{t} \geq \mathbf{0}, (\mathbf{t},\mathbf{h}) \leq H\}$, and we have $p_T = P_Y(\tilde{T}S_\mathbf{h})$.

Now we are in the conditions of Theorem 1. The boundary of the region $S_\mathbf{h}$ is regular, $B_Y \geq 0$ and $Y$ is continuous, since $X$ has these properties. The continuity of $X$ is guaranteed by the ss-property and by the condition (b). Indeed, using (2.2) and notation $E\xi^2 = \|\xi\|^2$, we have for $\mathbf{t}, \mathbf{s} \subset K$:

$$\|X(\mathbf{t}) - X(\mathbf{s})\| = \left\|\sum_{k=1}^d (X(t_1,...t_k s_{k+1},...s_d) - X(t_1,...t_{k-1} s_k...s_d))\right\|$$

$$\leq \sum_1^d \|X_i(t_i) - X_i(s_i)\| \leq \sqrt{c}\sum_1^d |t_i - s_i|^{\alpha_i} \leq C|\mathbf{t} - \mathbf{s}|^{\min \alpha_i}. \qquad (2.18)$$

The obtained estimate gives the continuity of $X$ in $K$ due to the Kolmogorov criterion [1].

By virtue of Theorem 1, the exponent $\theta_Y(S_\mathbf{h})$ exists and is finite, and therefore relation (2.17) can be written as follows

$$\theta_Y(S_\mathbf{h}) \geq \theta_X^+(K). \qquad (2.19)$$

Given the condition (d), we can state that $\theta_Y(S_\mathbf{h}) > 0$. Indeed, by (d), $B_{Y_i}(t) \leq c(1 \wedge t^{-(1+\varepsilon)})$; therefore $B_Y(\mathbf{t}) = \prod_i B_{Y_i}(t_i) \in L_1(R^d)$ and the desired estimate, $\theta_Y(S_\mathbf{h}) > 0$, follows from Theorem 1.

*Step* 2: the inverse inequality to (2.19).

Let $(\mathbf{H}_T, \|\cdot\|_{\mathbf{H}_T})$ be the Hilbert space with the reproducing kernel associated with $X(\mathbf{t})$ on $TK \subset R^d$. By Proposition 7, there exists $\varphi(\mathbf{t}) \in \mathbf{H}_T$ such that

$$\varphi(\mathbf{t}) \geq 1, \mathbf{t} \in U_T^C \quad \text{and} \quad \|\varphi\|_{\mathbf{H}_T} \leq c \ln^{(d-1)/2} T. \qquad (2.20)$$

(For $U_T^C$, see (2.15)). We have

$$P_X(TK) = P(X(\mathbf{t}) - \varphi_T(\mathbf{t}) \leq 1 - \varphi_T(\mathbf{t}), \mathbf{t} \in [0,T]^d)$$

$$\leq P(X(\mathbf{t}) - \varphi_T(\mathbf{t}) \leq 0, \mathbf{t} \in U_T^c) := p_T^{(\varphi)} \qquad (2.21)$$

To get rid of the trend $\varphi_T(\mathbf{t})$, we use the following inequality (see [3, 11]):

$$\left|\sqrt{-\ln p_T^{(\varphi)}} - \sqrt{-\ln p_T^{(0)}}\right| \leq \|\varphi_T\|_{\mathbf{H}_T}/\sqrt{2}. \qquad (2.22)$$

Combining this inequality with (2.21) we obtain

$$\sqrt{-\ln P_X(TK)} \geq \sqrt{-\ln p_T^{(\varphi)}} \geq \sqrt{-\ln P(X \leq 0, U_T^c)} - \|\varphi_T\|_{\mathbf{H}_T}/\sqrt{2}.$$

We divide both parts of this inequality into $\ln^{d/2} T$ and go to the limit on $T \to \infty$. Given that $\|\varphi\|_{\mathbf{H}_T}^2 = o(\ln^d T)$, we obtain

$$\theta_X^-(K) \geq \theta_Y(S_\mathbf{h}) = \theta_X(U_T^c, \ln^d T, 0). \qquad (2.23)$$

Combining (2.19) and (2.23), we have

$$\theta_Y(S_{\mathbf{h}}) \geq \theta_X^+(K) \geq \theta_X^-(K) \geq \theta_Y(S_{\mathbf{h}}),$$

i.e., the exponent $\theta_X(K)$ exists and is equal to $\theta_Y(S_{\mathbf{h}})$.

***Step 3***: the case $d = 2$. If $d = 2$, $S_{\mathbf{h}}$ is a right triangle with the legs $H/h_i, i = 1,2$. The rectangle $\Pi = (0, H/h_1) \times (0, H/h_2)$ can be divided diagonally into two figures corresponding to $S_{\mathbf{h}}$. They are compatible under the action of the group of movements $\Gamma$ described above. Therefore $S_{\mathbf{h}}$ is a parquet tile for any region with the regular boundary. Due to (2.5),

$$\theta_Y(S_{\mathbf{h}}) = \theta_Y(K)|S_{\mathbf{h}}| = \theta_Y(\sqrt{|S_{\mathbf{h}}|}K), \qquad |S_{\mathbf{h}}| = H^2/(2h_1 h_2).$$

## 3. Propositions: proof

**Proof of Proposition 2.** Let $Y(\mathbf{t})$ be a Gaussian random field with the covariance $B_Y \geq 0$ and homogeneous with respect to some group $\Gamma$ of movements of $R^d$; let $G_2$ be a region having the property of a parquet tile for the region $G_1$ (with respect to $\Gamma$). Then, for large T, we can find a set $\{\gamma_i, i = 1 - N_T\} \subset \Gamma$ such that

$$TG_1 \subseteq \cup_1^{N_T} \gamma_i G_2 := \hat{G}_{1T}, \qquad |TG_1| = N_T |G_2|(1 + o(1)).$$

Hence, for any $k > 0$, $P_Y(TkG_1) \geq P_Y(k\hat{G}_{1T})$. Given Slepian's lemma and the homogeneity of $Y(\mathbf{t})$, we have

$$P_Y(TkG_1) \geq P_Y(\cup_1^{N_T} \gamma_i kG_2) \geq \prod_{i=1}^{N_T} P_Y(\gamma_i kG_2) = [P_Y(kG_2)]^{N_T},$$

or, taking into account that $N_T|G_2| = |TG_1|(1 + o(1))$, we have

$$-\ln P_Y(TkG_1)/(Tk)^d \leq -|G_1|/|G_2|\ln P_Y(kG_2)/k^d (1 + o(1)).$$

Passing to infinity first over $T$, and then over $k$, we get

$$\theta_Y^+(G_1)/|G_1| \leq \theta_Y^-(G_2)/|G_2|.$$

Since the property of parquet tiles is mutual for $G_i, i = 1,2$, we shall also have

$$\theta_Y^+(G_2)/|G_2| \leq \theta_Y^-(G_1)/|G_1|.$$

Combining these inequalities, we get

$$\theta_Y^+(G_1)/|G_1| \leq \theta_Y^-(G_2)/|G_2| \leq \theta_Y^+(G_2)/|G_2| \leq \theta_Y^-(G_1)/|G_1|.$$

Hence $\theta_Y^+(G_i) \leq \theta_Y^-(G_i)$. This relation implies both the existence of exponents $\theta_Y(G_i), i = 1,2$ and the equality $\theta_Y(G_1)/|G_1| = \theta_Y(G_2)/|G_2|$.

**Proof of Proposition 6.**

We have to evaluate $EM_X(U_T)$ for $U_T = \{\mathbf{t} = (t_1,...t_d) : t_1^{h_1} \times ... \times t_d^{h_d} \leq 1\} \cap [0,T]^d$.

According to Dudley (see e.g.[1,Thm.1.3.3]), we have

$$EM_X(U_T) \leq c\int_0^\sigma \sqrt{\ln N(\varepsilon)}d\varepsilon, \qquad (3.1)$$

where

$$\sigma^2 = \max_{U_T} E[X(\mathbf{t})]^2 = \max_{U_T} \prod t_i^{h_i} = 1$$

and $N(\varepsilon)$ is the minimum number of open balls of radius $\varepsilon$ that are necessary to cover $U_T$. The balls are considered in the semi-metric $\rho(\mathbf{t},\mathbf{s}) = [E(X(\mathbf{t}) - X(\mathbf{s}))^2]^{1/2}$.

Similarly to (2.17),

$$\rho(\mathbf{t},\mathbf{s}) \leq \sum_{k=0}^{d-1} \rho(\mathbf{t}^{(k+1)}, \mathbf{t}^{(k)}), \quad \mathbf{t}^{(k)} = (t_1,...t_k, s_{k+1},...s_d).$$

Using (2.2) and assumption (b), we get

$$\rho(\mathbf{t},\mathbf{s}) \leq \sum_k t_1^{h_1}...t_{k-1}^{h_{k-1}} T^{h_k} s_{k+1}^{h_{k+1}}...s_d^{h_d} (E|X_k(t_k/T) - X_k(s_k/T)|^2)^{1/2}$$

$$\leq CT^H \|(\mathbf{t}-\mathbf{s})/T\|_\infty^{\alpha_-}, \qquad \alpha_- = \min \alpha_i, \; \|\mathbf{t}\|_\infty = \max|t_i|$$

Hence, $\rho(\mathbf{t},\mathbf{s}) \leq \varepsilon$ as soon as

$$\|\mathbf{t}-\mathbf{s}\|_\infty \leq (\varepsilon/C)^{1/\alpha_-} T^{1-H/\alpha_-} := \rho_\varepsilon \qquad (3.2)$$

To cover $U_T \subset TK$ with d-cubes of size $\rho_\varepsilon$, we will need cubes in the amount

$$N(\varepsilon) \leq T^d/\rho_\varepsilon^d = c(T^H/\varepsilon)^{d/\alpha_-}. \qquad (3.3)$$

Therefore, (3.3) is a rough estimate of the number of balls to cover $U_T$ in the $\rho(\mathbf{t},\mathbf{s})$ metric.
Hence

$$EM_X(U_T) \leq c_1 \int_0^1 \sqrt{\ln N(\varepsilon)} d\varepsilon = \tilde{c} T^H \int_0^{T^{-H}} \sqrt{\ln 1/x} dx$$

$$= \sqrt{\ln T}(C + o(1)) \qquad (3.4)$$

**Proof of Proposition 7.**

We need to find a function $\varphi(\mathbf{t}) \in \mathbf{H}_T$ with suitable properties (2.14) in the region $U_T^c = \{\mathbf{t} = (t_1,...t_d) : t_1^{h_1} \times ... \times t_d^{h_d} > 1\} \cap [0,T]^d$ ; $(\mathbf{H}_T, \|\cdot\|_{\mathbf{H}_T})$ is a Hilbert space with the reproducing kernel associated with $X(\mathbf{t})$ in $[0,T]^d$.

We define new coordinates: $\{\tilde{t}_i = h_i \ln t_i, i = 1-d\}$. Then the closure of $U_T^c$ looks as follows:

$$\tilde{U}_T^c = \{\tilde{\mathbf{t}} : (\tilde{\mathbf{t}}, \mathbf{e}) \geq 0, \; \tilde{\mathbf{t}} \leq \tilde{T}\mathbf{h}\} \subseteq \tilde{T}\prod_1^d \otimes [h_i - H, h_i], \qquad (3.5)$$

where $\tilde{T} = \ln T$, $\mathbf{e} = (1,...,1)$.

Let $Z = \{\mathbf{z}^{(i)}, i = 1-N_{\tilde{T}}\}$ be the set of integer points lying on the face $\pi_0 = \{\tilde{\mathbf{t}} \in \tilde{U}_T^c : (\tilde{\mathbf{t}}, \mathbf{e}) = 0\}$ of the simplex $\tilde{U}_T^c$. Obviously, the number of such points is $N_{\tilde{T}} \propto \tilde{T}^{d-1}$.

Let $\gamma(\tilde{\mathbf{t}}_0) = \{\tilde{\mathbf{t}} \geq \tilde{\mathbf{t}}_0\}$ be a cone with vertex $\tilde{\mathbf{t}}_0$. We're going to show that

$$\cup_{i=1}^{N_{\tilde{T}}} \gamma(\mathbf{z}^{(i)} - d\mathbf{e}) \supseteq \tilde{U}_T^c. \tag{3.6}$$

Let $\hat{\mathbf{t}}$ be the orthogonal projection $\tilde{\mathbf{t}} \in \tilde{U}_T^c$ on the hyperplane $(\tilde{\mathbf{t}}, \mathbf{e}) = \mathbf{0}$. Obviously, there is always an integer point $\mathbf{z}$ on $\pi_0$ such that

$$|\hat{t}_i - z_i| \leq 1, i = 2,...d \; ; \; |\hat{t}_1 - z_1| \leq d - 1. \tag{3.7}$$

But then $\hat{\mathbf{t}}$ and $\tilde{\mathbf{t}}$ belong to $\gamma(\mathbf{z} - d\mathbf{e})$. Indeed, $\tilde{\mathbf{t}} = \hat{\mathbf{t}} + (\tilde{\mathbf{t}}, \mathbf{e})\mathbf{e}/d$ and therefore

$$\tilde{t}_i - z_i + d_i = \hat{t}_i - z_i + d_i + (\tilde{\mathbf{t}}, \mathbf{e})\mathbf{e}/d \geq (\tilde{\mathbf{t}}, \mathbf{e})\mathbf{e}/d \geq 0.$$

Consider the random variable $\eta_T = A\sum_{\mathbf{z} \in Z} X(e^{(\mathbf{z}-d\mathbf{e})_i/h_i}, i=1-d)$. Taking into account that $(\mathbf{z},\mathbf{e}) = 0$, we have

$$E\eta_T^2 = A^2 e^{-2d^2} \sum_{Z \times Z} \prod_{i=1}^d B_{Y_i}((n_i - m_i)/h_i),$$

where $(\mathbf{n},\mathbf{m}) \in Z \times Z$ and $Y_i$ is the Lamperti transform of $X_i$.

Note that, under assumption (d), we have

$$\sum_{n=0}^{\infty} B_{Y_i}(n/h_i) < \infty, \tag{3.8}$$

because $B_{Y_i}(n/h_i) = B_{X_i}(1, e^{n/h_i})e^{-n_i} \leq n^{-1-\varepsilon}$, $n > n_0$.

To break ties, $(\mathbf{n},\mathbf{e}) = 0 = (\mathbf{m},\mathbf{e})$, between the components of the vectors $\mathbf{n}$ and $\mathbf{m}$, we exclude the variables $n_d$ and $m_d$ using the estimates: $0 \leq B_{Y_d}(n_d - m_d) \leq 1$. Given $B_{Y_i}(\cdot) \geq 0, i = 1-d$ and $(|n_i|,|m_i|, i = 1-d) \leq C\tilde{T}$ (see (3.5)), we have

$$E\eta_T^2 \leq A^2 e^{-2d^2} \prod_{i=1}^{d-1} \sum_{(|n_j|,|m_j|)<C\tilde{T}} B_{Y_i}((n_j - m_j)/h_i) \leq c\tilde{T}^{d-1} \tag{3.9}$$

Due to (3.8),

$$\sum_{(|n_j|,|m_j|)<C\tilde{T}} B_{Y_i}((n_j - m_j)/h_i) \leq c\tilde{T}.$$

Therefore $E\eta_T^2 \leq C\tilde{T}^{d-1}$.

Let us show that $\varphi(\mathbf{t}) = EX(\mathbf{t})\eta_T = EX(\mathbf{t})\hat{\eta}_T$ is the desired function, where $\hat{\eta}_T = E\{\eta_T|X(\mathbf{t}), t \in [0,T]^d\}$. Indeed, by definition $\varphi(\mathbf{t}) \in \mathbf{H}_T$; moreover,

$$\|\varphi\|_{\mathbf{H}_T} = E\hat{\eta}_T^2 \leq (E\eta_T^2)^{1/2} \leq c\tilde{T}^{(d-1)/2} = c(\ln^{(d-1)/2} T).$$

It remains to show that $\varphi(\mathbf{t}) \geq 1$, $\mathbf{t} \in U_T^c$.

Assume that $\tilde{\mathbf{t}} \in \tilde{U}_T^c$ belongs to a cone $\gamma(\mathbf{z} - d\mathbf{e})$. Then $\tilde{\mathbf{t}}$ can be represented as $\mathbf{z} - d\mathbf{e} + \mathbf{r}$, with $\mathbf{r} \phi 0$.

Since the covariance of $X$ is nonnegative

$$\varphi(\mathbf{t}) \geq A E X(e^{(z_i - d)/h_i}, i = 1 - d) \, X(e^{(z_i - d + r_i)/h_i}, i = 1 - d)$$

$$= Ae^{-2d^2} EX(1,...,1)) X(e^{r_1/h_1},...,e^{r_d/h_d}) = Ae^{-2d^2} \prod_{i=1}^{d} B_i(1, e^{r_i/h_i}) \ .$$

Therefore $\varphi(\mathbf{t}) \geq 1$, provided that $B_{X_i}(1, 1+x) \geq b_i$, $x > 0$ (assumption (e)) and $A = e^{2d^2} / \prod_{i=1}^{d} b_i$.

## 4. Exact values of $\theta_X$ : examples

We give a non-trivial example of $X = X_1 \otimes X_2$ with the exact value of the persistence exponent. The field is defined on the $Z_+ \times R_+$ space and has the covariance $B_X = (\delta_n^m + 1)/2 \cdot B(t,s)$. The process $X_1(n)$ is formally identical with discrete-time fractional Brownian motion of index h=0, whereas $X_2(t)$ is a Gaussian ss-process with $B(t,s) \geq 0$. We can consider $B_X$ as the covariance of a vector-valued process $[X_2^{(n)}(t) - X_2^{(0)}(t)]/\sqrt{2}$, where $\{X_2^{(n)}(t), 0 \leq n \leq T\}$ are independent copies of $X_2(t)$. Therefore

$$P(X(n,t) \leq 1, (n,t) \in (0.T)^2) = P(X_2^{(n)}(t) \leq X_2^0(t) + \sqrt{2}, 0 < n, t < T) := p_T. \tag{4.1}$$

The asymptotics of $p_T$ arise in the so-called pursuit problem [8]. In this context, the probability (4.1) means that during the period $(0,T)$ particles with dynamics $X_2^{(n)}(t)$, $X_2^{(n)}(0) = 0$, $n = 1,..., N = [T]$ will not be able to catch up with the particle moving according to $X_2^{(0)}(t) + \sqrt{2}$. Under some regular conditions on $B(t,s)$, [14],

$$\ln p_T = -\theta_X \ln^2 T(1 + o(1)) \ , \ \theta_X = [2\pi f(0)]^{-1} \ , \tag{4.2}$$

where $f(\lambda)$ is the spectral density of the stationary double of the process $X_2$. It follows that (4.2) is the persistence exponent in our example. In the case of the fractional Brownian motion, $X_2 = W_h$, the exponent is

$$\theta_X = \Gamma(1+h)/[2\Gamma(2h)\Gamma(1-h)] \geq 0.5h \wedge (1-h). \tag{4.3}$$

We show how (4.2) can be used to estimate $\theta_X$ from below for the field $X = X_{h_1} \otimes X_{h_2}$ with h-ss Gaussian components and nonnegative covariance. Let $Y = Y_{h_1} \otimes Y_{h_2}$ be the dual homogeneous field. Obviously,

$$P_Y(\tilde{T}K) \leq P(Y(i\Delta, s) \leq 0, (i\Delta, s) \in \tilde{T}K) \ , \ \tilde{T} = \ln T . \tag{4.4}$$

Assuming that $B_{Y_{h_1}}(t) \leq 1/2, t \geq \Delta$, we have

$$B_Y(i\Delta, s) \leq (\delta_0^i + 1)/2 \cdot B_{Y_{h_2}}(s) \tag{4.5}$$

Therefore, according to Slepian's lemma, the $Y(i\Delta, s)$ field in (4.4) can be replaced by a field with the covariance $(\delta_0^i + 1)/2 \cdot B_{Y_{h_2}}(s)$. But then, due to (4.2), we get $\theta_Y(K) \geq [2\pi f_{Y_2}(0)]^{-1}$

Finally, we have

$$\theta_X(K) = \theta_Y(K) \cdot 2(h_a/h_g)^2 \geq [2\pi f_{Y_2}(0)]^{-1} \cdot 2(h_a/h_g)^2 . \tag{4.6}$$

For the (fractional) Brownian sheet, the relations (4.3, 4.6) give

$$\theta_X(K) \geq h \wedge (1-h), \quad X = W_h \otimes W_h . \tag{4.7}$$

For comparison, we note the result for the fractional Brownian motion $W_h^{(d)}$ with d-dimensional time for which $E|W_h^{(d)}(\mathbf{t}) - W_h^{(d)}(\mathbf{s})|^2 = |\mathbf{t}-\mathbf{s}|^{2h}$, $W_h^{(d)}(\mathbf{0}) = 0$, $0 < h < 1$. In this case the order of decay of the persistence log-probability, $\psi_T = \ln T$, is independent of $d$, and the exponent is known explicitly: $\theta_{W_h^{(d)}}([0,T] \times [-T,T]^{d-1}, \ln T, 1) = d - h$, [13].

The problem of choosing $\psi_T$ is well studied now for Gaussian stationary processes [7]. Examples of exact values of the exponents $\theta_X$ are still rare. The general state of the persistence problem for processes is presented in [4, 5].